\documentclass{amsart}

\usepackage[utf8]{inputenc}

\usepackage{amssymb,tikz}
\usepackage{amsthm,dsfont}
\usepackage{amsmath, a4wide}
\usepackage{amsbsy}
\usepackage[all]{xy}
\usepackage{bm}
\usepackage{hyperref}
\usepackage{color}
\usepackage{subfigure}
\usepackage{tikz}

\newtheorem{thm}{Theorem}[section]
\newtheorem{cor}[thm]{Corollary}

\theoremstyle{definition}
\newtheorem{rem}[thm]{Remark}

\newcommand{\NN}{\mathbb{N}}

\newcommand{\RR}{\mathbb{R}}

\DeclareMathOperator{\so}{sopfr}
\DeclareMathOperator{\sod}{sopf}

\author{Dimitris Vartziotis$^{1,2}$}
\address{$^1$ TWT GmbH Science \& Innovation, Department for Mathematical Research, Ernsthaldenstr. 17, 70565 Stuttgart, Germany }

\author{Aristos Tzavellas$^2$}
\address{$^2$ NIKI Ltd. Digital Engineering, Research Center, 205 Ethnikis Antistasis Street,
45500 Katsika, Ioannina, Greece$$ $$}

\email{dimitris.vartziotis@twt-gmbh.de }

\title[On sums of prime factors]{On sums of prime factors}

\begin{document}
%\author{Dimitris Vartziotis}
%\author{Aristos Tzavellas}
\begin{abstract} We study the arithmetic function $\so(n)$ (OEIS A001414) which gives the sum of prime factors (with repetition) of a number $n$. In particular we obtain the asymptotic formula
$$ \sum_{n \leq x} \so(n) \sim \frac{\pi^2}{12} \frac{x^2}{\log x},$$
which holds as well for the function $\sod(n)$ (OEIS A008472) that just gives the sum of distinct prime factors of $n$.
This asymptotic formula was already stated by R. Jakimcyuk \cite{rj12} which was brought to our attention after the completion of the first version of this manuscript.\\
%{\bf Last Update:} 28/3/2017
\end{abstract}
\maketitle
\pagestyle{myheadings}
\thispagestyle{empty}
\markboth{DIMITRIS VARTZIOTIS AND ARISTOS TZAVELLAS}{On sums of prime factors}

%%%%%%%%%%%%%%%%%%%%%%%%%%%%%%%%%%%%%%%%%%%%%%%%%%%%%%%%%%
\section{Introduction}
\label{sec1}
The number and distribution of primes $p$ less than a given number $x \in \RR$ is a classical problem in number theory. Let $\pi(x) = \sum_{p \leq x} 1,$ $x\in \RR$, be the prime counting function and write $F(x) \sim G(x)$ if $\lim_{x \rightarrow \infty} F(x)/G(x)=1$. By the prime number theorem we have that
$$ \pi(x) \sim \frac{x}{\log x}. $$ 
Using the prime number theorem it can easily be shown that 
\begin{equation} \label{sum1}
P(x)=\sum_{p \leq x} p \sim \frac{1}{2} \frac{x^2}{\log x};
\end{equation}
see f.e. \cite[Section 2.7]{bach}. The aim of this note is to extend the sum on the left hand side of \eqref{sum1} to all integers in the following way. For $n \in \NN$, $n > 1$, we write
$$ n = \prod_{i=1}^r p_i^{\alpha_i} = p_1^{\alpha_1} \cdot p_2^{\alpha_2} \cdot \ldots \cdot p_r^{\alpha_r}, $$
for the unique prime factorization of $n$, in which the $p_i$ are the different prime factors of $n$ and the $\alpha_i$ are the corresponding multiplicities. We define the arithmetic function
$$ \so(n) = 
\begin{cases} 
\sum_{i=1}^r \alpha_i p_i = \alpha_1 p_1 + \alpha_2 p_2 + \ldots + \alpha_r p_r &\mbox{for } n > 1, \\ 
0 & \mbox{for } n = 1.
\end{cases} 
$$
Clearly, $n$ is prime if and only if $\so(n)=n$. Moreover, $\so(n)$ is completely additive due to the uniqueness of the prime factorization of every integer $n$.
The function $\so(n)$ fluctuates wildly attaining sharp local maxima whenever $n$ is prime and small values whenever $n$ has many small prime factors. As an example, assume that $n$ is a Mersenne prime, i.e., $n=2^p - 1$ for a prime $p$, then $\so(n) = 2^p - 1$ and $\so(n+1) = \so(2^p) = 2p$. As a second example, consider the numbers $10^9 + 7$, $10^9 + 8$, $10^9 + 9$. The first and the last of these numbers are primes and we obtain:
\begin{align*}
\so (10^9 + 7) &= 10^9 + 7 \\
\so (10^9 + 8) &= 3 \cdot 2 + 2 \cdot 3 + 1 \cdot 7 + 2 \cdot 109 + 1 \cdot 167 = 404\\
\so (10^9 + 9) &= 10^9 + 9 
\end{align*}
It is thus interesting to study 
$$ B(x) = \sum_{n \leq x} \so(n) $$
and compare it with the growth of $P(x)$. The value of $B(x)$ is obviously larger than $P(x)$ for all $x \geq 4$. The main result of this note is the following asymptotic formula for $B(x)$.

\begin{thm} \label{thm1}
$$B(x) \sim \frac{\pi^2}{12} \frac{x^2}{\log x} $$
\end{thm}
We can relate this formula in a straightforward way to the function $P(x)$.
\begin{cor}
$$B(x) \sim \frac{\pi^2}{6} P(x) $$
\end{cor}
The interesting observation is that $\pi^2/6 = 1.664 \ldots$ is considerably larger than $1$, however, the order of the growth of the two sums is the same.

%%%%%%%%%%%%%%%%%%%%%%%%%%%%%%%%%%%%%%%%%%%%%%%%%%%%%%%%%%
\section{Proof of Theorem \ref{thm1}}
\label{sec2}

For a real number $x$ we write $[x]$ for the integer part of $x$. 
We start from the definition $B(x) = \sum_{n \leq x} \so(n)$ and rewrite $B(x)$ as a sum over all primes $p \leq x$. We count the contribution of each $p$ to $B(x)$ in a systematic way as follows. First, we note that exactly $[x/p]$ integers $n \leq x$ have $p$ as a prime factor, from these integers exactly $[x/p^2]$ are also divisible by $p^2$, and from these exactly $[x/p^3]$ are divisible by $p^3$ as well. Continuing in this manner we see that the contribution of each $p$ to $B(x)$ is $p ( [\frac{x}{p}] + [\frac{x}{p^2}] + [\frac{x}{p^3}] + \ldots )$. Therefore,
\begin{equation} \label{sumMain}
B(x) = \sum_{p \leq x} p \left ( \left[\frac{x}{p} \right] + \left[\frac{x}{p^2} \right] + \left[\frac{x}{p^3} \right] + \ldots \right ).
\end{equation}
Now,
\begin{align*} 
\sum_{p \leq x} p \left[\frac{x}{p^2} \right] &\leq \sum_{p \leq x} p \frac{x}{p^2} = x \sum_{p\leq x} \frac{1}{p} \\
&= x (\log \log x + \mathcal{O}(1)) = x \log \log x + \mathcal{O}(x),
\end{align*}
by a classical result of Mertens \cite{mertens} and 
\begin{align*}
\sum_{p \leq x} p \left ( \left[\frac{x}{p^3} \right] + \ldots \right ) & \leq \sum_{p \leq x} p \left ( \frac{x}{p^3} + \ldots  \right ) =  x \sum_{p \leq x} \left( \frac{1}{p^2} + \ldots \right ) \\
&= x \sum_{p \leq x}  \frac{1}{p^2} \frac{1}{1 - \frac{1}{p}} = x \sum_{p \leq x}  \frac{1}{p(p-1)} = \mathcal{O}(x) 
\end{align*}
since $\sum_{p \leq x}  \frac{1}{p(p-1)}$ converges. Therefore, we are left with the estimation of 
\begin{equation} \label{sum2}
\sum_{p \leq x} p \left[ \frac{x}{p} \right ].
\end{equation}
\begin{rem}
Note that the trivial estimate $\sum_{p \leq x} p \left[ \frac{x}{p} \right ] \leq x \pi(x)$ satisfies $x \pi(x) \sim x^2/ \log x$.
\end{rem}
\begin{rem}
From our argument it is clear that \eqref{sum2} can be used to rewrite $\sum_{n\leq x} \sod(n)$, in which $\sod(n)$ gives the sum of distinct prime factors of $n$. Thus, the order of growth is the same as for $B(x)$.
\end{rem}
Recalling that $P(x) = \sum_{p \leq x} p$, we see that the contribution of the primes contained in the interval $(\frac{x}{2},x]$ to the sum \eqref{sum2} is
$$ \sum_{\frac{x}{2} < p \leq x} p = P(x) - P\left(\frac{x}{2} \right). $$
The contribution of the primes in  $(\frac{x}{3}, \frac{x}{2}]$ is
$$2 \sum_{\frac{x}{3} < p \leq \frac{x}{2}} p = 2 \left( P\left(\frac{x}{2} \right) - P\left(\frac{x}{3}\right) \right), $$
and the contribution of the primes in  $(\frac{x}{4}, \frac{x}{3}]$ is
$$3 \sum_{\frac{x}{4} < p \leq \frac{x}{3}} p = 3 \left( P\left(\frac{x}{3} \right) - P\left(\frac{x}{4}\right) \right).$$
Thus, 
\begin{align*}
\sum_{p \leq x} p \left[ \frac{x}{p} \right ] &= P(x) - P \left (\frac{x}{2} \right ) + 2 \left ( P \left (\frac{x}{2} \right ) - P \left (\frac{x}{3} \right )\right ) + 3 \left( P \left (\frac{x}{3} \right ) - P \left (\frac{x}{4} \right ) \right) + \ldots \\
&= P(x) + P \left (\frac{x}{2} \right ) + P \left (\frac{x}{3} \right ) + \ldots = \sum_{1 \leq n \leq \frac{x}{2} } P \left( \frac{x}{n} \right).
\end{align*} 

We would like to use this formula to obtain an asymptotic formula for \eqref{sum2}. Therefore, we need an estimate for the error term in \eqref{sum1}.
Starting from the prime number theorem in the form 
\begin{equation} \label{primeCount}
\pi(x) = \frac{x}{\log x} + \mathcal{O}\left (  \frac{x}{\log^2 x} \right ),
\end{equation}
we use partial summation to obtain
\begin{equation} \label{sum3}
P(x) = \frac{x^2}{ 2 \log x} + \mathcal{O} \left ( \frac{x^2}{\log^2 x} \right).
\end{equation}
Next, we use this formula to rewrite our expression as follows.
\begin{align}
\notag \sum_{1 \leq n \leq \frac{x}{2} } P \left( \frac{x}{n} \right) &= \sum_{1 \leq n \leq \frac{x}{2} } \left ( \frac{1}{2} \left( \frac{x}{n} \right)^2 \frac{1}{\log \frac{x}{n}}  + \mathcal{O}\left( \left( \frac{x}{n} \right)^2 \frac{1}{\log^2 \frac{x}{n}} \right) \right ) \\
\label{sum4} &= \frac{x^2}{2} \sum_{1 \leq n \leq \frac{x}{2} } \frac{1}{n^2 \log \frac{x}{n} } + x^2 \mathcal{O} \left( \sum_{1 \leq n \leq \frac{x}{2} } \frac{1}{n^2 \log^2 \frac{x}{n}} \right).
\end{align}
To estimate the first sum, we use again partial summation. We set
$$ u(n) = \frac{1}{n^2} \ \ \ \ \text{ and } \ \ \ \ f(n) = \left( \log \frac{x}{n} \right)^{-1}. $$
Then,
\begin{align*} 
U(x) &= \sum_{n \leq x} u(n) = \sum_{n \leq x} \frac{1}{n^2} = \frac{\pi^2}{6} + \mathcal{O}\left( \frac{1}{x} \right) \\
f'(n) &= (-1) \left( \log \frac{x}{n} \right)^{-2} \frac{n}{x} \left ( - \frac{x}{n^2} \right) = \frac{1}{n \log^2 \frac{x}{n} }.
\end{align*}
Consequently,
\begin{align*}
\sum_{1 \leq n \leq \frac{x}{2} } \frac{1}{n^2 \log \frac{x}{n} } &= \sum_{n \leq \frac{x}{2} } u(n) f(n) = U\left( \frac{x}{2} \right) f\left( \frac{x}{2} \right) - \int_1^{\frac{x}{2} } U(t) f'(t) \ dt \\
&= \left( \frac{\pi^2}{6} + \mathcal{O}\left( \frac{1}{x} \right) \right) \frac{1}{\log 2} - \int_1^{\frac{x}{2} } \left( \frac{\pi^2}{6} + \mathcal{O}\left( \frac{1}{t} \right) \right) \frac{1}{t \log^2 \frac{x}{t} } \ dt \\
&= \frac{\pi^2}{6 \log 2}  + \mathcal{O}\left( \frac{1}{x} \right) - \frac{\pi^2}{6}  \int_1^{\frac{x}{2} } \frac{dt}{t \log^2 \frac{x}{t} } + \mathcal{O} \left( \int_1^{\frac{x}{2} } \frac{dt}{t^2 \log^2 \frac{x}{t} } \right ) \\
&=\frac{\pi^2}{6} \frac{1}{\log x} + \mathcal{O} \left ( \frac{1}{\log^2 x} \right ),
 \end{align*}
since
\begin{align*}
\int_1^{\frac{x}{2} } \frac{dt}{t \log^2 \frac{x}{t} } &\overset{u=\frac{x}{t}} = \int_x^{2 } \frac{u}{x} \frac{1}{\log^2 u} \left( - \frac{x}{u^2} \right) \ du \\
&= \int_{2}^{x} \frac{du}{u \log^2 u} = \left . - \frac{1}{\log u} \right \vert_2^x = \frac{1}{\log 2} - \frac{1}{\log x},
\end{align*}
and similarly
\begin{align*}
\int_1^{\frac{x}{2} } \frac{dt}{t^2 \log^2 \frac{x}{t} } = \mathcal{O}\left (  \frac{1}{\log^2 x} \right).
\end{align*}

Setting $f(t) = (\log \frac{x}{t})^{-2}$ we can handle the sum in the $\mathcal{O}$ term of \eqref{sum4} in a similar fashion and obtain
\begin{equation}\label{sum5}
\sum_{1 \leq n \leq \frac{x}{2} } \frac{1}{n^2 \log^2 \frac{x}{n} } = \frac{\pi^2}{6} \frac{1}{\log^2 x} + \mathcal{O}\left( \frac{1}{\log^3 x} \right) = \mathcal{O}\left( \frac{1}{\log^2 x} \right).
\end{equation}
Putting the pieces together and using \eqref{sum4} we get
$$ \sum_{1 \leq n \leq \frac{x}{2} } P \left( \frac{x}{n} \right) = \frac{\pi^2}{12} \frac{x^2}{\log x} + \mathcal{O}\left( \frac{x^2}{\log^2 x} \right). $$
Consequently, we obtain the right asymptotic order of $B(x)$ via \eqref{sumMain}
\begin{align} 
\notag B(x) &= \frac{\pi^2}{12} \frac{x^2}{\log x} + \mathcal{O} \left( \frac{x^2}{\log^2 x} \right)  + \mathcal{O}(x \log \log x) + \mathcal{O}(x) \\
\label{final} &= \frac{\pi^2}{12} \frac{x^2}{\log x} + \mathcal{O} \left( \frac{x^2}{\log^2 x} \right).
\end{align}

\begin{rem}
Note that even under the assumption of the Riemann Hypothesis (RH) the error term coming from the prime number theorem will be dominant in \eqref{final}. Assuming the RH we can replace the error term in \eqref{primeCount} by the slightly sharper term; see f.e. \cite{granville} for details.
\end{rem}

\section*{Acknowledgement}

We would like to thank Florian Pausinger from Technische Universität München for enhancing the presentation of the paper.

%%%%%%%%%%%%%%%%%%%%%%%%%%%%%%%%%%%%%%%%%%%%%%%%%%%%%%%%%%%%%%%%%
%%%%%%%%%%%%%%		Bibliography
%%%%%%%%%%%%%%%%%%%%%%%%%%%%%%%%%%%%%%%%%%%%%%%%%%%%%%%%%%%%%%%%%

\end{document}